\documentclass[10pt,a4,twoside,rm]{article}
\usepackage[left=4.5cm,top=4cm,right=4.5cm,bottom=4cm]{geometry}
\usepackage{amsmath, amsthm, amssymb}
\usepackage{graphicx}
\usepackage{amssymb}
\usepackage{fancyhdr}
\usepackage{titlesec}
\usepackage{youngtab}
\usepackage{color}
\usepackage{hyperref}
\usepackage{pstricks}
\usepackage{pst-all}

\newenvironment{dem}{\noindent \textit{Proof:} }{\quad \hfill $\square$}

\newcommand{\hm}{\operatorname{Hom}}
\newcommand{\en}{\operatorname{End}}

\newcommand{\rad}{\operatorname{rad}}
\newtheorem{teo}{Theorem}[section]
\newtheorem{lem}[teo]{Lemma}

\newtheorem{pro}[teo]{Proposition}
\newtheorem{defi}[teo]{Definition}
\newtheorem{cor}[teo]{Corollary}

\newtheorem{rem}[teo]{Remark}
\newtheorem{con}[teo]{Conjecture}
\newtheorem{cla}[teo]{Claim}
\newcommand{\K}{\mathbb K}
\newcommand{\R}{ \mathbb R}

\newcommand{\Z}{ \mathbb Z}

\newcommand{\co}[1]{\text{coef}_{1^{\otimes}}(#1)}
\newcommand{\bol}[1]{ \boldsymbol{#1}}
\newcommand{\kl}[1]{ \underline{H}_{#1}}
\newcommand{\bs}[1]{ BS(\sub{#1})}
\newcommand{\bsa}[1]{ BS(#1)}
\newcommand{\cell}[1]{ \Delta_{\sub{w}}(#1)}
\newcommand{\A}{ A_{\sub{w}} }
\newcommand{\sub}[1]{ \underline{#1} }
\newcommand{\gd}{\dim_{v}}
\newcommand{\simpleset}[1]{\Lambda_0(\sub{ #1})}

\newcommand{\simple}[1]{ D_{\sub{w}}(#1)  }
\newcommand{\gdn}[2]{d_{\sub{w}}(#1,#2)}

\newcommand{\lr}[1]{ \langle #1 \rangle}
\numberwithin{equation}{section}

\titleformat{\section}[hang]{\sc}{\thesection.}{0.5cm}{\filcenter}
\titleformat{\subsection}[hang]{\bf}{\thesubsection.}{0.2cm}{\filright}

\begin{document}

\title{ \normalsize {\bf \sc Graded Cellularity and the Monotonicity Conjecture}}
\author{\sc  David Plaza }
\date{ }
\maketitle

\begin{center}
Universidad de Chile,\\
Facultad de Ciencias,\\
Santiago, Chile
\end{center}

\begin{center}
Email: davidricardoplaza@gmail.com
\end{center}

\begin{abstract}
\noindent \textsc{Abstract.} The graded cellularity of Libedinsky double leaves, which form a basis for the endomorphism ring of the Bott-Samelson-Soergel bimodules, allows us to view the Kazhdan-Lusztig polynomials as graded decomposition numbers. Using this interpretation, we can provide a  proof of the monotonicity conjecture for any Coxeter system.
\end{abstract}

\section{Introduction}

In their seminal paper \cite{kl}, Kazhdan and Lusztig defined, for each Coxeter system $(W, S)$, a family of polynomials with integer coefficients indexed by pairs of elements of $W$. These polynomials are now known as the Kazhdan-Lusztig (KL) polynomials. We will denote them by $P_{x, w}(q) \in \Z[q]$, for all $x, w \in W$. Applications of the KL-polynomials have been found in the representation theory of semisimple algebraic groups, the  topology of Schubert varieties, the theory of Verma modules, the Bernstein-Gelfand-Gelfand (BGG) category $\mathcal{O}$, etc. (see, e.g., \cite{br} and references therein).

\medskip
Aside from the importance of the KL-polynomials in the above-mentioned subjects, there are purely combinatorial reasons to study these polynomials. Perhaps the major reason is the longstanding Kazhdan-Lusztig positivity conjecture \cite{kl}, which states that $P_{x, w}(q) \in \mathbb{N}[q]$, for all Coxeter groups $W$ and all $x, w \in W$. In 2013, Elias and  Williamson \cite{ew1} gave a proof of this conjecture by proving a stronger result known as Soergel's conjecture.

\medskip
For each Coxeter group $W$, Soergel constructed a category of graded $R$-bimodules (where $R$ is a polynomial ring with coefficients in $\R$) known as the category of Soergel bimodules, which we will denote by $\mathbb{S} \text{Bim}$. He proved that (up to degree shift) $W$ parameterizes the set of all indecomposable objects in $\mathbb{S} \text{Bim}$. For $w \in W$, let us denote by $B_w$ the corresponding indecomposable object. Soergel proved in \cite{so} that $\mathbb{S} \text{Bim}$ is a categorification of the Hecke algebra $\mathcal{H}$ of $W$. This means that there is an algebra isomorphism
\begin{equation}  \label{eq intro one}
\epsilon: \mathcal{H} \rightarrow [\mathbb{S}  \text{Bim}],
\end{equation}
where $[\mathbb{S} \text{Bim}]$ denotes the split Grothendieck group of $\mathbb{S} \text{Bim} $. Soergel proposed the following
conjecture, which came to be known as Soergel's conjecture:
\begin{equation}    \label{eq intro two}
  \epsilon (\underline{H}_{w}) = [B_{w}],
\end{equation}
where $\{ \underline{H}_{w}  \}_{w \in W}$ is the Kazhdan-Lusztig basis of $\mathcal{H}$. Assuming this conjecture, Soergel showed that the coefficients of the Kazhdan-Lusztig polynomials of $W$ are given by the dimensions of certain Hom spaces in $\mathbb{S} \text{Bim}$. It follows that these coefficients are non-negative, i.e., he proved that (\ref{eq intro two}) implies the positivity conjecture.

\medskip
More than mere positivity is conceivable for coefficients of KL-polynomials. In effect, Braden and MacPherson \cite[Corollary 3.7]{bm} proved a monotonicity property for these coefficients, when $W$ is a finite or an affine Weyl group, by using the notion of sheaves on moment graphs. Concretely, if $W$ is a finite or an affine Weyl group and $u, v, w \in W$ such that $u \leq v \leq w$, then

\begin{equation}  \label{eq intro}
P_{u, w}(q) - P_{v, w}(q) \in \mathbb{N}[q],
\end{equation}
where $\leq $ denotes the usual Bruhat order on $W$. In other words, if we fix the second index of a KL-polynomial, and if the first one decreases in Bruhat order, all coefficients in the polynomial weakly increase in value. It is natural to conjecture that (\ref{eq intro}) holds for arbitrary Coxeter groups. In the literature, the latter conjecture is referred to as the \emph{Monotonicity Conjecture} for KL-polynomials.

\medskip
Although moment graphs originated in geometry, Fiebig observed in \cite{fi} that it is possible to define them axiomatically. In particular, he associated a moment graph with any Coxeter system. Using this point of view, Fiebig gave an alternative description of $\mathbb{S} \text{Bim}$. It then seems reasonable to think that the arguments given by Braden and MacPherson in the finite and affine Weyl group settings can be generalized in order to prove the Monotonicity Conjecture for arbitrary Coxeter groups. However, it is not clear to the author how to confirm that the results obtained by Elias and Williamson regarding Soergel bimodules (Soergel's conjecture) translate to the moment graph setting. It should be mentioned that although the above is probably known to experts, it has not been formally documented anywhere, to the best of the author's knowledge.

\medskip
In this paper, we provide a proof of the monotonicity conjecture completely contained in the language of Soergel bimodules. Furthermore, our proof does not refer to Fiebig's theory of sheaves on moment graphs. Therefore, the arguments used in our proof are different from the ones used by Braden and MacPherson in the finite and affine Weyl group cases. Let us briefly explain our approach to the Monotonicity Conjecture. For each reduced expression $\sub{w}$ of an element $w \in W$, one can explicitly define a Soergel bimodule $\bs{w}$, called the Bott-Samelson bimodule. The endomorphism ring of a Bott-Samelson bimodule, $\text{End}(\bs{w})$, has a natural structure of free right $R$-algebra. Libedinsky constructed in \cite{li1} an $R$-basis for these spaces that he called \emph{light leaves basis}. He generalized his construction in \cite{li} to obtain another basis that he called the \emph{double leaves basis}. The latter is more useful than the light leaves basis for our purposes because of its symmetry properties. In particular, Elias and Williamson proved in \cite{ew} that the double leaves basis is a cellular basis for $\text{End}(\bs{w})$ in the sense of Graham and Lehrer \cite{gl}.

\medskip
Let $R^{+}$ be the ideal of $R$ generated by homogeneous elements of nonzero degree. We have $\R \cong R/R^+$. Therefore, we can reduce  $\text{End}(\bs{w})$ modulo $R^{+}$ to obtain an $\R$-algebra. The resulting algebra is equipped with a natural $\Z$-grading. The double leaves basis behaves satisfactorily with respect to reduction modulo $R^{+}$ and cellularity. Concretely, the image of the double leaves basis is a graded cellular basis of $\text{End}(\bs{w})\otimes_{R} \R$ in the sense of Hu and Mathas \cite{hm}. The existence of a graded cellular basis allows us to define graded cell modules and graded simple modules, as well as graded decomposition numbers. We then prove using Soergel's conjecture that the KL-polynomials (suitably normalized) can be interpreted as graded decomposition numbers. Finally, we construct certain injective homomorphisms between cell modules that allow us to embed one cell module into the other. This embedding implies a monotonicity property for the respective graded decomposition numbers that is equivalent to the Monotonicity Conjecture according to the aforementioned interpretation of the KL-polynomials as graded decomposition numbers.

\medskip
The layout of this article is as follows. In Section 2, we recall a few useful results of the theory of graded cellular algebras. In Section 3, we define Hecke algebras and the category of Soergel bimodules, and conclude this section by recalling Libedinsky's construction of the double leaves basis. We establish the graded cellularity of double leaves basis in Section 4. Using graded cellularity, we can view the KL-polynomials as graded decomposition numbers. Finally, in Section 5, we show how to embed a cell module into another. We then use this embedding to conclude our proof of the Monotonicity Conjecture.

\section{Graded cellular algebras}  \label{section gca}
In this section, we briefly recall the theory of graded cellular algebras. Graded cellular algebras were defined by Hu and Mathas in \cite{hm}, following and extending the ideas of Graham and Lehrer \cite{gl}. A clear exposition of this theory (in the ungraded setting) can be found in \cite{ma}.

\medskip
Let $\K$ be a field. A graded $\K$-vector space $M$ is a $\K$-vector space that has a direct sum decomposition $M = \bigoplus_{k\in \Z} M_{k}$. If $M$ is a graded $\K$-vector space and $k \in \Z$, we denote by $M\langle k \rangle$ the graded $\K$-vector space obtained from $M$ by shifting the grading on $M$, i.e., $M\langle k \rangle_{i} = M_{i - k}$, for all $i \in \Z$. Given a Laurent polynomial $f = \sum_{i \in \Z} a_i v^{i} \in \mathbb{N}[v,v^{-1}]$ and a graded vector space $M$, we set

\begin{equation*}
fM = \bigoplus_{i\in \Z}  M\lr{i}^{\bigoplus{a_i}}.
\end{equation*}

A graded $\K$-algebra $A$ is a $\K$-algebra with a direct sum decomposition $A = \bigoplus_{k \in \Z}A_{i}$ as a $\K$-vector space such that $A_{i}A_{j} \subset A_{i + j}$, for all $i, j \in \Z$. A graded right $A$-module $M$ is a graded $\K$-vector space that is an $A$-module in the usual (ungraded) sense, such that $A_{i}M_{j}\subset M_{i + j}$, for all $i, j \in \Z$. Let $v$ be an indeterminate over $\Z$. If $M = \bigoplus_{k \in \Z}, M_{k}$ is a graded finite-dimensional $\K$-vector space. We define its graded dimension as the following Laurent polynomial:

\begin{equation}
\gd M = \sum_{k\in \Z} \dim_{\K} (M_{k})v^k
\end{equation}

We now define the concept of graded cellular algebra. This definition is provided in \cite[Definition 2.1]{hm}
\begin{defi} \rm   \label{definition gca}
Let $A$ be a graded finite-dimensional $\K$-algebra. A graded cell datum is an ordered quadruple $(\Lambda, T, C ,\deg)$, where $(\Lambda, \geq)$ is a poset, $T(\lambda) $ is a finite set for $\lambda \in \Lambda$, and $C$ and $\deg $ are two functions defined as follows:

$$ C:\coprod_{\lambda \in \Lambda} T(\lambda) \times T(\lambda) \rightarrow A \text{, } (\mathfrak{s},\mathfrak{t})\rightarrow c_{\mathfrak{st}}^{\lambda}; \qquad \qquad \deg:\coprod_{\lambda \in \Lambda} T(\lambda) \rightarrow \Z     $$
 such that $C$ is injective and:
\begin{description}
\item[\rm(a)] $\mathcal{C} = \{ c_{\mathfrak{st}}^{\lambda} \mbox{ } | \mbox{ } \mathfrak{s}, \mathfrak{t} \in T(\lambda), \mbox{ } \lambda \in \Lambda  \}$ is a basis of $A$.
\item[\rm(b)]  The $\K$-linear map $\ast: A\rightarrow A$ determined by $(c_{\mathfrak{st}}^{\lambda} )^{\ast} = c_{\mathfrak{ts}}$ is an anti-automorphism of $A$.
\item[\rm(c)] For all $a \in A $, $\lambda \in \Lambda$, and $\mathfrak{s},\mathfrak{t}\in T(\lambda)$, there exist scalars $r_{\mathfrak{tv}}(a)\in \K$ that do not depend on $\mathfrak{s}$, such that
 \begin{equation} \label{Aaction}
 c_{\mathfrak{st}}^{\lambda}a\equiv \sum_{\mathfrak{v}\in T(\lambda)} r_{\mathfrak{tv}}(a) c_{\mathfrak{sv}}^{\lambda} \mod A^{> \lambda}
 \end{equation}
 where $A^{> \lambda} $ is the vector subspace of $A$ spanned by $\{ c_{\mathfrak{ab}}^{\mu} \mbox{ }| \mbox{ } \mathfrak{a,b} \in T(\mu), \mu >\lambda \}$.
 \item[\rm(d)] Each $c_{\mathfrak{st}}^{\lambda}$ is a homogeneous element of degree $\deg(\mathfrak{s}) + \deg(\mathfrak{t})$.
\end{description}
 A graded cellular algebra is a graded algebra with a graded cell datum. The set $\mathcal{C}$ is a graded cellular basis of $A$.
 \end{defi}

\begin{rem} \label{remark gca whitout grading} \rm
Ignoring the grading on $A$, the degree function, and axiom (d) in the above definition, we can recover the original definition of cellular algebras by Graham and Lehrer \cite{gl}. In this case, we say that $A$ is a cellular algebra with a cellular basis and a cell datum.
\end{rem}

Let $A$ be a graded cellular algebra with graded cellular basis $\mathcal{C}$, as in the above definition. For each $\lambda \in \Lambda$, we define the \emph{graded cell module}, $\Delta(\lambda)$, as the graded right $A$-module
$$ \Delta(\lambda) = \bigoplus_{k\in \Z} \Delta(\lambda)_{k}   $$
where $\Delta(\lambda)_{k}$ is a $\K$-vector space with basis $\{ c_{\mathfrak{t}}^{\lambda} \mbox{ }|\mbox{ } \mathfrak{t} \in T(\lambda) \text{ and } \deg(\mathfrak{t}) = k  \}$, and where the $A$-action on $\Delta(\lambda)$ is determined by the scalars that appear in (\ref{Aaction}), i.e.,
\begin{equation}
c_{\mathfrak{t}}^{\lambda} a = \sum_{\mathfrak{v}\in T(\lambda)} r_{\mathfrak{tv}}(a)c_{\mathfrak{v}}^{\lambda}
\end{equation}

Suppose that $\lambda \in \Lambda$. Then, it follows from Definition \ref{definition gca} that there is a bilinear form $\langle \mbox{ } , \mbox{ } \rangle$ on $\Delta(\lambda)$ which is determined by
\begin{equation}
c_{\mathfrak{as}}^{\lambda} c_{\mathfrak{tb}}^{\lambda} \equiv  \lr{c_{\mathfrak{s}}^{\lambda}, c_{\mathfrak{t}}^{\lambda}} c_{\mathfrak{ab}}^{\lambda}  \mod A^{>\lambda}
\end{equation}

For each $\lambda \in \Lambda$, $\lr{\mbox{ } , \mbox{ }}$ is a symmetric and associative bilinear form on $\Delta(\lambda)$, where radical $\rad(\Delta(\lambda))$ is a graded $A$-submodule (see \cite[Section 2]{hm}). Therefore, the quotient $D(\lambda):= \Delta(\lambda)/\rad(\Delta(\lambda))$ is a graded right $A$-module. Furthermore, if $D(\lambda)\neq 0$, then $D(\lambda)$ is a simple graded right $A$-module. Define
$$ \Lambda_0 = \{ \lambda \in \Lambda \mbox{ }| \mbox{ } D(\lambda)\neq 0      \} .$$

The following theorem gives a classification of the simple graded A-modules for a graded cellular algebra $A$. This result is due to Hu and Mathas \cite[Theorem 2.10]{hm}, and is a graded version of \cite[Theorem 3.4]{gl}.
\begin{teo} \label{cellular theory}
Let $A$ be a graded cellular algebra, with a cell datum as in Definition \ref{definition gca}. Then, the set $\{ D(\lambda)\langle k \rangle \mbox{ } | \mbox{ } \lambda \in \Lambda_0 \text{ and } k\in \Z \}$ is a complete set of pairwise non-isomorphic graded simple right $A$-modules.
\end{teo}

Let $\Delta$ and $D$ be graded right $A$-modules. If $D$ is simple, we denote by $[\Delta: D\lr{k}]$ the multiplicity of the graded simple module $D\lr{k}$ as a graded composition factor in a graded composition series of $\Delta$, for all $k \in \Z$. We then define the graded decomposition number, $d(\Delta, D)$, as the Laurent polynomial:
\begin{equation}
d(\Delta, D) = \sum_{k\in Z} [\Delta : D\lr{k}]v^{k}
\end{equation}
In particular, if $\Delta = \Delta(\lambda)$ and $D = D(\mu)$, for some $\lambda \in \Lambda$ and $\mu \in \Lambda_0 $, we denote
\begin{equation}
d(\Delta, D) = d(\lambda, \mu)
\end{equation}

We end this section by relating the graded representation theory of algebras $A$ and $eAe$, where $A$ is a graded (not necessarily cellular) algebra and $e\in A$ is homogeneous idempotent. For each right $A$-module $V$, the subspace $Ve$ of $V$ has a natural structure of a right $eAe$-module.

\begin{teo}  \label{idempotent algebra}
Let $A$ be a graded algebra. Let $e \in A$ be an homogeneous idempotent (and, therefore, of degree zero). We then have:
\begin{description}
\item[\rm (a)] If $V$ is a simple graded right $A$-module and $Ve\neq 0$, $Ve$ is a simple graded right $eAe$-module. Furthermore, all the simple right $eAe$-modules can be obtained in this manner.
\item[\rm (b)] Let $V$ and $D$ be graded right $A$-modules. If $D$ is simple and $De\neq 0$, then
  \begin{equation}  \label{gdn equal}
  d(V, D) = d(Ve, De)
  \end{equation}
 where the left (resp. right) side of \ref{gdn equal} corresponds to the graded decomposition number for $A$-modules (resp. $eAe$-modules).
\item[\rm (c)] If $e$ is primitive idempotent, there is a unique (up to degree shift) simple graded right $eAe$-module.
\item[\rm (d)] If $e$ is primitive idempotent, $D$ is the unique simple $eAe$-module, and $\gd D = 1$, then $d(V, D) = \gd V$.
\end{description}
 \begin{dem}  \rm
 The first and second claims can be found in \cite[Appendix A1]{do} and \cite[Theorem 2.4]{pl}, respectively. For the last claim, note first that since $e$ is primitive idempotent, $P =: eA$ is a principal indecomposable right $A$-module and $P/\text{rad}\mbox{ } P$ is a simple $A$-module. It is a well-known fact that for an $A$-module $M$, $P/ \text{rad} \mbox{ }P$ is a composition factor of $M$ if and only if $Me\neq 0$ (see, for example, \cite[Proposition A14]{ma}). Suppose now that $M$ and $N$ are simple $A$-modules such that $Me\neq 0$ and $Ne\neq 0$. Then, $P/\text{rad}\mbox{ } P$ is a composition factor of $M$ and $N$; however, since $M$ and $N$ are simple modules, we have $M\cong N$. Then, (c) follows from (a) above. Part (d) is clear.
 \end{dem}

\end{teo}

\section{Libedinsky double leaves}   \label{section dlb}
In this section, we introduce, for an arbitrary Coxeter system $(W, S)$, its corresponding Hecke algebra, and its corresponding category of Soergel bimodules. We end this section by introducing the double leaves basis. This is a basis for morphism spaces between  Bott-Samelson bimodules. Double
leaves are the combinatorial tool that we use to prove the monotonicity conjecture in Section \ref{section monotony}. This basis admits a convenient diagrammatic description. For the sake of brevity, we have omitted the diagrammatic approach in this paper. However, the diagrams allowed several calculations that helped us understand the problem. We refer the reader interested in the diagrammatic approach to \cite[Part 3]{ew}.
\subsection{Hecke algebras and KL-polynomials}
Let $(W, S)$ be a \emph{Coxeter system}. That is, $W$ is a group with generators $s\in S$ and relations
\begin{equation}  \label{relations}
(st)^{m_{st}} = e     \qquad  \qquad    \text{ for all } s, t\in S
\end{equation}
where $e \in W$ is the identity, $m_{st}\in \{1, 2, \ldots, \infty\}$ satisfies: $m_{st} = 1$ if and only if $s = t$, and $m_{st} = m_{ts}$ for all $s, t \in S$. When $m_{st} = \infty $, relation (\ref{relations}) is omitted. We use the \emph{underlined} letter $\sub{w} = (s_{1}, \ldots, s_k) $, $s_i \in S$ to denote a finite sequence of elements in $S$. We will call \emph{expressions} to these sequences. If we consider an expression $\sub{w} = (s_{1}, \ldots , s_k)$, the corresponding Roman letter, $w$, will denote its product in $W$, i.e., $w = s_1 \ldots s_{k}$.  We make this distinction between $w$ and $\sub{w}$ because a few concepts defined in this paper and used throughout rely heavily on the considered expression for $w$ rather than on $w$ itself. We will often write $\sub{w} = s_{1} \ldots s_{k}$, where the underlined letter reminds us that the entire sequence, and not merely $w$, is important. The group $W$ is equipped with a length function $l:W \rightarrow \mathbb{N}$ and an order, called the \emph{Bruhat order}, which is denoted by $\geq$ (see, e.g.,  \cite[Chapter 1]{bb}). The \emph{length} of an  expression $\sub{w} = s_1 \ldots s_{k}$ is $k$. We say that an expression is reduced if $l(\sub{w}) = l(w)$.

\begin{defi} \label{Hecke algebra}
Let $\mathcal{A} = \mathbb{Z}[v, v^{-1}]$ be the ring of the Laurent polynomials in $v$. The Hecke algebra $\mathcal{H} = \mathcal{H}(W, S)$ is the $\mathcal{A}$-algebra that is associative and unital with generators $\{  H_{s}  | s\in S \}$ and relations
\begin{equation}
H_{s}^{2} = (v^{-1}-v)H_s + 1
\end{equation}
\begin{equation}
\underbrace{H_{s}H_{t}H_{s} \ldots}_{m_{st} \text{ terms }} = \underbrace{H_{t}H_{s}H_{t} \ldots}_{m_{st} \text{ terms }}
\end{equation}
\end{defi}
If $\sub{w} = s_{1} \ldots s_{k}$ is a reduced expression for $w\in W$, we define $H_{w} := H_{s_{1}} \ldots H_{s_{k}}$. It is well-known that $H_{w}$ does not depend on the choice of the reduced expression $\sub{w}$. The set $\{ H_{w} |w\in W \}$ is a basis for $\mathcal{H}$ as an $\mathcal{A}$-module. There is a unique ring involution $-: \mathcal{H}\rightarrow \mathcal{H}$ determined by $\overline{v} = v^{-1}$ and $\overline{H_{w}} = H_{w^{-1}}^{-1}$, for all $w\in W$.

\begin{teo} \label{kl theorem}
There exists a unique basis $\{ \kl{w} \mbox{ } | \mbox{ } w\in W  \}$ for $\mathcal{H}$ as a $\mathcal{A}$-module such that $\kl{w}$ is invariant under $-$ and
\begin{equation}
\kl{w} = \sum_{x\leq w} h_{x, w}H_{x}
\end{equation}
with $h_{x, w} \in v\mathbb{Z}[v]$ if $x\neq w$ and $h_{w, w} = 1$.
\end{teo}

The set $\{ \kl{w} \mbox{ } | \mbox{ } w\in W\}$ is called the \emph{Kazhdan-Lusztig basis} of $\mathcal{H}$ and the polynomials $h_{x, w}$ are called the \emph{Kazhdan-Lusztig polynomials}.

\begin{rem} \rm \label{variables}
The reader should note that in this paper, we follow the normalization given by Soergel in \cite{so} rather than the original normalization by Kazhdan and Lusztig in \cite{kl}. Therefore, we have $q = v^{-2}$, and the original Kazhdan-Lusztig polynomials $P_{x, w}(q)\in \mathbb{Z}[q]$ can be recovered from our Kazhdan-Lusztig polynomials $h_{x, w}(v)\in \mathbb{Z}[v]$ by the formula
\begin{equation}
h_{x,w}(v) = v^{l(w) - l(x)}P_{x, w}(v^{-2})
\end{equation}
\end{rem}

\subsection{The category of Soergel bimodules} \label{the category of bimodule of soergel}

Let us fix once and for all a \emph{reflection-faithful} representation $V$ of $W$ over $\R$. In \cite{so1}, Soergel constructed such a representation for arbitrary Coxeter groups. Let $R$ be the $\R$-algebra of  regular functions on $V$. We can grade this algebra by setting $R = \bigoplus_{i\in \mathbb{Z}}R_i$, with $R_2 = V^{\ast}$. Let $R^{+}$ be the ideal of $R$ generated for all elements of positive degree. Of course, $R / R^+ \cong \R$. We will often consider $\R$ as an $R$-module via this isomorphism. There is a natural action of $W$ on $R$ induced by the action of $W$ on $V$. For $s \in S$, let $R^{s}$ be the subring of $R$ fixed by $s$. Then, we define the graded $(R, R)$-bimodule
\begin{equation}
B_{s} = R \otimes_{R^{s}} R (1),
\end{equation}
where for a graded $(R, R)$-bimodule $B$ and every $k \in \Z$, we denote by $B(k)$ the graded $(R, R)$-bimodule defined by the formula
\begin{equation}  \label{shift bimodules}
B(k)_{i} = B_{k + i}
\end{equation}

For the expression $\sub{w} = s_{1} \ldots s_{s_{k}}$, we denote by $\bs{w}$ the $(R, R)$-bimodule defined by
\begin{equation}
\bs{w} = B_{s_{1}} \otimes_{R} B_{s_{2}} \otimes_R \ldots \otimes_{R} B_{s_{k}}
\end{equation}
Bimodules of the type $\bs{w}$ will be called \emph{Bott-Samelson bimodules}. We introduce the convention that $BS(\emptyset) = R$. From now on, we denote the tensor product of $(R, R)$-bimodules, $\otimes_{R}$, simply juxtaposition. Thus, $\bs{w}$ becomes $B_{s_1}B_{s_2} \ldots B_{s_{k}}$. We then have the following isomorphism of $(R, R)$-bimodules
\begin{equation}
B_{s_1}B_{s_2} \ldots B_{s_{k}} \cong R \otimes_{R^{s_1}} R \otimes_{R^{s_2}} \otimes \ldots \otimes_{R^{s_{k}}} R  (k)
\end{equation}
Therefore, we can write an element of this module as a sum of terms given by $k + 1$ polynomials in $R$, one in each slot separated by the tensors. Let $x_{s} \in V^\ast$ be an equation of the hyperplane fixed by $s \in S$. Then, for all $s \in S$, we define the Demazure operator, $\partial_{s}: R (2) \rightarrow R^s$, as a morphism of graded $R^{s}$-modules given by
\begin{equation}
\partial_{s}(f) = \frac{f - s \cdot f}{2x_{s}}
\end{equation}
It is not difficult to prove that $\partial_{s}(f)$ and  $P_s(f) = f - x_s\partial_s(f)$ are $s$-invariant. Since $f = P_s(f) + x_s\partial_s(f)$, $R$ is free as a graded right $R^{s}$-module with basis $\{1, x_s\}$, i.e., we have a decomposition $R\simeq R^{s}\oplus x_s R^s$. Using this decomposition, we can prove that $BS(s)$ is a right free $R$-module with basis $\{ x_s \otimes 1, 1\otimes 1  \}$.  Let $\underline{w} = s_1\ldots s_k$ be an expression. Going once more through the above decomposition of $R$, we can see that $\bs{w}$ is a right free $R$-module of rank $2^{k}$, with basis

\begin{equation} \label{base of bott-samelson}
 \{  x_{s_1}^{e_1}\otimes x_{s_2}^{e_2} \otimes \ldots \otimes x_{s_k}^{e_k}\otimes 1 \mbox{ } | \mbox{ } e_{i} = 0, 1  \}.
\end{equation}
We now define the category of Soergel bimodules that categorifies the Hecke algebra, as will be made
precise in Theorem \ref{soergel cat theo}.

\begin{defi}  \rm
The category of \emph{Soergel bimodules}, $\mathbb{S}\text{Bim} $, is the category of $\mathbb{Z}$-graded $(R, R)$-bimodules whose objects are grading shifts and direct sums of direct summands of Bott-Samelson bimodules. The morphisms are all degree-preserving bimodule homomorphisms. For $B, B' \in \mathbb{S}\text{Bim}$, we denote by $\hm(B, B')$ the corresponding set of morphisms. We also define
\begin{equation}
\hm^{\Z}(B, B') = \bigoplus_{k\in \Z} \hm(B(k), B')
\end{equation}
An element $f\in \hm(B(k), B') $ is called a homogeneous map of degree $k$.
\end{defi}

Let $[\mathbb{S}\text{Bim}]$ be the split Grothendieck group of the category $\mathbb{S}\text{Bim}$. That is, $[\mathbb{S}\text{Bim}]$ is the abelian group generated by the symbols $[B]$ for all objects $B \in \mathbb{S}\text{Bim}$, subject to the relation $[B] = [B'] + [B'']$ whenever we have $B\cong B'\oplus B''$ in $\mathbb{S}\text{Bim} $. The following theorem is known as Soergel's categorification theorem and relates $\mathbb{S}\text{Bim}$ to $\mathcal{H}$.

\begin{teo}    \label{soergel cat theo}
For each $w \in W$, there exists a unique (up to isomorphism) indecomposable bimodule $B_w$ that occurs as a direct summand of $\bs{w}$ for any reduced expression $\underline{w}$ of $w$, and $B_w$ does not appear in any $\bs{x}$ for a word $\underline{x}$ shorter than $\underline{w}$. Furthermore, there is a unique $\mathcal{A}$-algebra isomorphism
\begin{equation}
\epsilon: \mathcal{H} \longrightarrow [\mathbb{S}\text{Bim}]
\end{equation}
such that $\epsilon (v) = R( 1 )$ and $\epsilon(\kl{s}) = [B_s]$ for all $s\in S$.
\end{teo}
In order to explicate the inverse of $\epsilon$, known as Soergel's character map, we need to introduce standard bimodules. Given $x \in W$ we define the \emph{standard bimodule} $R_{x}$ as the $(R, R)$-bimodule, such that $R_{x}\cong R$ as a left $R$-module and the right action on $R_{x}$ is the right multiplication on $R$ deformed by the action of $x$ on $R$, i.e.,
\begin{equation}
 r\cdot r' := rx(r') \text{ for } r\in R_x \text{ and } r'\in R.
\end{equation}

\begin{teo} \label{inverse categorification}
The categorification $\epsilon: \mathcal{H} \rightarrow [\mathbb{S}Bim]$ admits an inverse, $\eta:[\mathbb{S}Bim] \rightarrow \mathcal{H} $, given by the formula
\begin{equation}
\eta([B]) = \sum_{x\in W}  \gd (\hm^{\Z}(B, R_x)\otimes_{R}\R)   H_{x}
\end{equation}
\end{teo}

We end this subsection by introducing Soergel's conjecture. For historical reasons, we call this a \emph{conjecture} even though it was proven in 2013 \cite{ew1}.
\begin{con}  \label{soergel conjecture}
Let $W$ be a Coxeter group. For all $w \in W$, we have
 \begin{equation}  \label{soergel conjecture equation}
 \epsilon (\underline{H}_{w}) = B_{w}.
\end{equation}
\end{con}

\begin{rem} \rm
The Kazhdan-Lusztig positivity conjecture immediately follows from Soergel's conjecture by applying Soergel's character map to (\ref{soergel conjecture equation}).
\end{rem}

\subsection{Double leaves basis}

Let $\sub{w}$ and $\sub{v}$ be two (not necessarily reduced) expressions. In this section, we recall the construction of the double leaves basis (DLB), a basis of the space $\hm^{\Z}(\bs{w}, \bs{v} )$, defined by Libedinsky in \cite{li}. The DLB is, in some sense, an improvement over the light leaves basis, another basis for $\hm^{\Z}(\bs{w}, \bs{v} )$ defined in \cite{li1}. In the remainder of this paper, we will work with the DLB rather than the light leaves basis because as we will see in Section \ref{section cellularity}, DLB is a (graded) cellular basis whereas the light leaves basis is not. We use the cellularity of the DLB to establish the monotonicity conjecture for Kazhdan-Lusztig polynomials.

\medskip
To introduce the DLB, we begin by defining three morphisms between the Bott-Samelson bimodules. The first is the multiplication morphism, $m_{s}$, which is a degree $1$ morphism determined by the formula:
\begin{equation}
\begin{array}{r}
  m_{s}: \bsa{s} = R\otimes_{R^{s}} R (1)\rightarrow R \\
     p\otimes q \rightarrow pq
\end{array}
\end{equation}
The second morphism is a unique (up to multiplication by a nonzero scalar) $-1$ degree morphism, $j_{s}$, determined by the formula
\begin{equation}
\begin{array}{r}
  j_{s}: \bsa{ss} = R\otimes_{R^{s}} R \otimes_{R^{s}} R (2) \rightarrow \bsa{s} \\
     1\otimes p \otimes 1 \rightarrow \partial_{s}(p) \otimes 1
\end{array}
\end{equation}

For $s, r\in S$, consider the bimodule
$$ X_{sr} := BS(srs\ldots)  $$
with the product having $m_{sr}$ terms. We then define $f_{sr}$ as the unique degree zero morphism from $X_{sr}$ to $X_{rs}$ sending $1\otimes 1 \otimes \ldots \otimes 1$ to $1\otimes 1\otimes \ldots \otimes 1$. We denote by $\mathbb{I}$ the identity on the endomorphism ring of a Bott-Samelson bimodule. Each time we use the symbol $\mathbb{I}$, the relevant Bott-Samelson bimodule will be cleared from the context. For each expression $\sub{w} = s_1 \ldots s_n \in S^{n}$, we inductively define a perfect binary directed tree, denoted by $\mathbb{T}_{\sub{w}}$, with nodes colored by Bott-Samelson bimodules and edges colored by morphisms from parent nodes to child nodes. At depth 1, the tree looks as in Figure \ref{uno}.
\psset{unit=0.5}

\begin{figure}[h]
\centering{
\begin{pspicture}(0,0)(6,4)
\rput(3,3){\rnode{P1}{$\bsa{s_1\ldots s_n}$}}
\rput(-1,-2){\rnode{a1}{$\bsa{s_2\ldots s_{n}}$}}
\rput(7,-2){\rnode{a2}{$\bsa{s_1\ldots s_n}$}}
\rput(-1,1){$m_{s_1} \otimes \mathbb{I}^{n-1}$}
\rput(6,1){$\mathbb{I}$}
\ncline[nodesep=0.15cm]{->}{P1}{a1}
\ncline[nodesep=0.15cm]{->}{P1}{a2}
\end{pspicture} }
\vspace{1cm}
\caption{Level one of $\mathbb{T}_{\sub{w}}$}
\label{uno}
\end{figure}

Now, let $1 < k \leq n $ and assume that we have constructed the tree to level $k - 1$. Let $ \sub{u} = t_{1} \ldots t_{i} \in S^{i}$ be a node $N$ of depth $k - 1$ colored by the bimodule $\bsa{t_1 \ldots t_i} \bsa{s_{k}\ldots s_{n}}$. We then have two possibilities:
\begin{description}
\item[a) $l(t_1 \ldots t_is_{k}) > l(t_1 \ldots t_i) $.] In this case the child nodes and edges of $N$ are constructed as shown in Figure \ref{dos}.

\begin{figure}[h]
\centering{
 \begin{pspicture}(0,0)(6,4)
\rput(3,3){\rnode{P1}{$\bsa{t_1 \ldots t_i} \bsa{s_{k}\ldots s_{n}}$}}
\rput(-3,-2){\rnode{a1}{$\bsa{t_1 \ldots t_i} \bsa{s_{k+1}\ldots s_{n}}$}}
\rput(9,-2){\rnode{a2}{$\bsa{t_1 \ldots t_i} \bsa{s_{k}\ldots s_{n}}$}}
\rput(-3,1){$\mathbb{I}^{i} \otimes m_{s_k} \otimes \mathbb{I}^{n-k}$}
\rput(7,1){$\mathbb{I}$}
\ncline[nodesep=0.15cm]{->}{P1}{a1}
\ncline[nodesep=0.15cm]{->}{P1}{a2}
\end{pspicture}  } \vspace{1cm}
\caption{Level $k$ of $\mathbb{T}_{\sub{w}}$}
\label{dos}
\end{figure}

\item[b) $l(t_1 \ldots t_i s_{k}) < l(t_1 \ldots t_i) $.] In this case, it is a well-known fact for Coxeter groups that there exists a sequence of braid moves that converts $\sub{u} = t_1 \ldots t_i $ into $\sub{u'}=t'_1 \ldots t'_{i - 1} s_{k}$. Of course, there are several ways to do this. However, we can fix a particular sequence of braid moves and construct a morphism $\bsa{\sub{u}} \rightarrow \bsa{\sub{u'}}$ by replacing each braid move in the sequence by its respective morphism of type $f_{sr}$. We denote this morphism by $F(\sub{u}, \sub{u'}, s_{k})$. The child nodes of $N$ are then colored by the two Bott-Samelson bimodules located at the bottom of Figure \ref{tres}, and the child edges are colored by morphisms obtained by composing the dashed arrows in Figure \ref{tres}.

\begin{figure}[h] \label{tres}
\centering{
 \begin{pspicture}(-6,0)(6,14)
\rput(0,14){\rnode{P1}{$\bsa{t_1 \ldots t_i} \bsa{s_{k}\ldots s_{n}}$}}
\rput(0,10){\rnode{P2}{$\bsa{t'_1 \ldots t'_{i-1}s_k} \bsa{s_{k}\ldots s_{n}}$}}
\rput(0,6){\rnode{P3}{$\bsa{t'_1 \ldots t'_{i-1}s_k\ldots s_{n}}$}}
\rput(-6,0){\rnode{a1}{$\bsa{t'_1 \ldots t'_{i-1}} \bsa{s_{k+1}\ldots s_{n}}$}}
\rput(6,0){\rnode{a2}{$\bsa{t'_1 \ldots t'_{i-1}} \bsa{s_{k}\ldots s_{n}}$}}
\ncline[linestyle=dashed,nodesep=0.15cm]{->}{P1}{P2}
\ncline[linestyle=dashed,nodesep=0.15cm]{->}{P2}{P3}
\ncline[linestyle=dashed,nodesep=0.15cm]{->}{P3}{a1}
\ncline[linestyle=dashed,nodesep=0.15cm]{->}{P3}{a2}
\ncarc[arcangle=60,nodesep=0.15cm]{->}{P1}{a2}
\ncarc[arcangle=-60,nodesep=0.15cm]{->}{P1}{a1}
\rput(2.2,12){$\scriptstyle F(\sub{u},\sub{u'},s_k)\otimes \mathbb{I}^{n-k+1}$}
\rput(2.2,8){$\scriptstyle \mathbb{I}^{i-1}\otimes j_{s_k} \otimes \mathbb{I}^{n-k}$}
\rput(-4.3,4){$\scriptstyle \mathbb{I}^{i-1}\otimes m_{s_k} \otimes \mathbb{I}^{n-k}$}
\rput(3,4){$\scriptstyle \mathbb{I}$}
\end{pspicture} } \vspace{1cm}
\caption{Level $k$ of $\mathbb{T}_{\sub{w}}$}
\label{tres}
\end{figure}

\end{description}

By construction, each leaf of the tree $\mathbb{T}_{\sub{w}}$ is colored by a Bott-Samelson bimodule $\bs{u}$, where the expression $\sub{u}$ is reduced. Note that it is possible for two leaves to be colored by Bott-Samelson bimodules $\bs{u}$ and $\bs{u'}$, where $\sub{u}$ and $\sub{u'}$ are two reduced expressions for the same element in $u \in W$. To avoid this ambiguity, we realize the following choices:
\begin{enumerate}
  \item Fix, once and for all, a reduced expression $\sub{u}$ for all $u \in W$.
  \item For all $u \in W$ and all reduced expression $\sub{u'}$ of $u$, choose a sequence of braid moves that converts $\sub{u'}$ into $\sub{u}$, where $\sub{u}$ is the fixed reduced expression selected in the previous step.
  \item Finally, replace each braid move in the sequence selected in the previous step by its corresponding morphism of type $f_{sr}$ to obtain a morphism from $\bs{u'}$ to $\bs{u}$, denoted by $F(\sub{u'}, \sub{u})$.
\end{enumerate}

We now complete the construction of $\mathbb{T}_{\sub{w}}$ by composing each of the lower Bott-Samelson bimodules with its corresponding morphism  $F(\sub{u'}, \sub{u})$. This procedure avoids the ambiguity in coloring the leaves. That is, if two leaves are colored by $\bs{u}$ and $\bs{u'}$, where $\sub{u}$ and $\sub{u'}$ are two reduced expressions for the same element $u \in W$, $\sub{u}$ and $\sub{u'}$ are the same expressions.

\medskip
By composing the corresponding arrows, we can consider each leaf in $\mathbb{T}_{\sub{w}}$ colored by $\bs{u}$ as a morphism from $\bs{w}$ to $\bs{u}$, where $\sub{u}$ is the reduced expression for $u \in W$ fixed above. Let $\mathbb{L}_{\sub{w}}(u)$ be the set of all leaves colored by $\sub{u}$. As mentioned before, we will consider the set $\mathbb{L}_{\sub{w}}(u)$ as a subset of $\hm^{\Z}(\bs{w}, \bs{u})$. Note that every leaf is a homogeneous morphism, since it was constructed as a composition of homogeneous morphisms. In fact, the degree of each leaf can be computed as $+1$ for each occurrence of a morphism of type $m_s$ and $-1$ for each occurrence of a morphism of type $j_s$.

\begin{rem}  \label{dependence of the choice} \rm
The set $\mathbb{L}_{\sub{w}}(u)$ is not uniquely determined because it relies heavily on the choices realized along the way. Thus, when we refer to it, one must understand that we are implicitly assuming that we have fixed a particular choice for each of the non-canonical steps in the construction of $\mathbb{T}_{\sub{w}}$. For example, if $\sub{w}$ is a reduced expression for $w \in W$, there is exactly one leaf in $\mathbb{L}_{\sub{w}}(w)$. This leaf can be chosen as any morphism of the type $F(\sub{w}, \sub{w'})$, where $\sub{w'}$ is any reduced expression of $w$.  However, we choose identity in this case for the sake of simplicity. We will henceforth use this choice throughout the paper without reference to it.
\end{rem}

\medskip
In order to introduce the DLB, we need to define an adjoint leaf for each leaf. To do this, we must first define an adjoint morphism for each of $m_{s}, j_{s}$, and $f_{sr}$. The adjoint morphism of $m_{s}$ is
\begin{equation}  \label{epsilon}
    \begin{array}{rl}
          \epsilon_{s}: R  \rightarrow & \bsa{s} \\
          1  \rightarrow & x_{s}\otimes 1 + 1 \otimes x_{s}
        \end{array}
\end{equation}

For $j_{s}$, the corresponding adjoint morphism is
\begin{equation}
   \begin{array}{rl}
         p_{s}: \bsa{s}  \rightarrow & \bsa{ss} \\
          1\otimes 1  \rightarrow &  1\otimes 1 \otimes 1
        \end{array}
\end{equation}

Finally, for $f_{sr}$, the adjoint morphism is $f_{rs}$. Each adjoint morphism is homogeneous and has the same degree as its corresponding morphism. For each leaf  $l:\bs{w} \rightarrow \bs{u}$ in $ \mathbb{T}_{\sub{w}}$, we can thus define an adjoint leaf $l^{a}:\bs{u} \rightarrow\bs{w}$ as the morphism obtained by replacing each morphism of type $m_s, j_{s}$ and $f_{sr}$ by its corresponding adjoint. We thus obtain an inverted tree, $\mathbb{T}_{w}^a$, with the same nodes as $\mathbb{T}_{w}$ but with the arrows pointing in the opposite direction.

\medskip
For $f \in \hm^{\Z}(\bs{w},\bs{u} )$ and $g\in \hm^{\Z}(\bs{x}, \bs{y})$ we define
\begin{equation}
f\cdot g =\left\{
            \begin{array}{ll}
             f \circ g  , & \text{ if } \sub{w} = \sub{y}; \\
             \emptyset , & \text{ if } \sub{w} \neq \sub{y}.
            \end{array}
          \right.
\end{equation}

For an expression $\sub{w}$, we denote by $\mathbb{L}_{\sub{w}}$ the set of all leaves in $\mathbb{T}_{w}$. We are now in a position to define the main object of interest in this paper.

\begin{teo} \cite[Theorem 3.2]{li}
For all expressions $\sub{w}$ and $ \sub{v}$, the set $\mathbb{L}_{\sub{v}}^a  \cdot \mathbb{L}_{\sub{w}}$ is a basis as right $R$-module of the space $\hm^{\Z}(\bs{w}, \bs{v})$. We call this set the \emph{double leaves basis}.
\end{teo}

In order to prove the linear independence of $\mathbb{L}_{\sub{v}}^a  \cdot \mathbb{L}_{\sub{w}}$, Libedinsky \cite{li} introduced an order on the set $\mathbb{L}_{\sub{v}}^a  \cdot \mathbb{L}_{\sub{w}}$ and applied a classical triangularity argument. This order is defined by indexing each leaf by two sequences of zeros and ones, which we denote by $\bol{i} = (i_1, \ldots, i_{n})$ and $\bol{j} = (j_1, \ldots , j_n)$. Let us recall this assignment, since it will be important for our purposes.

\medskip
Let $\sub{w} = s_1 \ldots s_n$ be an expression of length $n$. Recall that $l$ was constructed inductively in $n$ steps. For all $1 \leq k \leq n$,  we set $i_{k} = 1$ if a morphism of type $m_{s}$ appears in the $k$-th step of the construction of $l$; otherwise, we set $i_{k} = 0$. In a similar manner, we set $j_{k} = 1$ if a morphism of type $j_{s}$ appears in the $k$-th step of the construction of $l$, and otherwise set $j_{k} = 0$. Note that each leaf $l$ is completely determined by these two sequences and the expression $\sub{w}$. Thus, we denote $l = f_{\bol{i}}^{\bol{j}}$. For $s \in S$, let us denote $x_{s}^{0} = 1$ and $x_{s}^1 = x_{s}$. If $\bol{j} = (j_1, \ldots, j_n)$ is a sequence of zeros and ones, we set

\begin{equation}
x^{\bol{j}} = x_{s_1}^{j_1} \otimes x_{s_2}^{j_2} \otimes \ldots \otimes x_{s_{n}}^{j_n} \otimes 1 \in \bs{s}
\end{equation}

In particular, if $\bol{j} = (0, \ldots , 0)$, we denote $x^{\bol{j}}$ by $1^{\otimes}$. The indexation of each leaf by pairs of binary sequences is compatible with the lexicographic order in the sense of the following lemma \cite[Section 5.5]{li}:

\begin{lem}   \label{lemma lex}
Denote by $\geq$ the lexicographic order on $\{0, 1\}^{n}$. Then,
$$ f_{\bol{i}}^{\bol{j}}(x^{\bol{j}'}) = \left\{
     \begin{array}{rl}
       1^{\otimes}, & \text{ if } \bol{j} = \bol{j}' \text{;}\\
       0, & \text{ if } \bol{j} > \bol{j}'.
     \end{array}
   \right.
$$
\end{lem}

We end this section by introducing a basis for $\hm^{\Z}(\bs{w}, R_x)$, for all reduced expressions $\sub{w}$ of $w \in W$ and $ x \in W$. We recall that $R_x$ denotes the standard bimodule defined immediately following Theorem \ref{inverse categorification}. We first need to introduce a new morphism. For all $s \in S$, consider the $(R, R)$-bimodule morphism  $\beta_s: BS(s) \rightarrow R_s$ determined by $\beta_s (p \otimes q) = ps(q)$, for all $p, q \in R$. For all $x, y \in W$, we have $R_xR_y \cong R_{xy}$. Therefore, we can also define a morphism that we denote by $\beta_x: \bs{x} \rightarrow R_{x}$, for all $x \in W$. Let us define the set
\begin{equation}  \label{standard basis equation}
\mathbb{L}_{\sub{w}}^{\beta}(x) = \{  \beta_x \circ l \mbox{ }| \mbox{ } l \in \mathbb{L}_{\sub{w}}(x)  \} \subset \hm^{\Z} (\bs{w}, R_x).
\end{equation}
Following Libedinsky, we will call $\mathbb{L}_{\sub{w}}^{\beta}(x)$ the standard leaves basis. This name is justified by the following lemma (see\cite[Proposition 6.1]{li}).

\begin{lem}  \label{standard basis}
Let $\sub{w}$ be an expression and let $x \in W$. Then, $\mathbb{L}_{\sub{w}}^{\beta}(x)$ is an $R$-basis of $\hm^{\Z} (\bs{w}, R_x)$ as a right $R$-module.
\end{lem}

\begin{cor} \label{corollary standard basis}
Let $\sub{w}$ be a reduced expression for $w$ and let $x \in W$. Then, $\hm^{\Z} (\bs{w}, R_x) \neq 0$ if and only if $x \leq w$.
\end{cor}
\begin{dem}
The result is a direct consequence of Lemma \ref{standard basis} once we note that $\mathbb{L}_{\sub{w}}(x)\neq \emptyset$ if and only if $x\leq w$.
\end{dem}

\section{KL-polynomials as graded decomposition numbers}   \label{section cellularity}
In this section, we interpret KL-polynomials as graded decomposition numbers. We first need to establish the graded cellularity of the double leaves basis. For the rest of this section, we fix a reduced expression $\underline{w}$ for an element $w \in W $. In \cite[Proposition 6.22]{ew}, Elias and Williamson observed that $\en^{\Z}(\bs{w})$ is a cellular $R$-algebra with the double leaves basis as cellular basis. Let us specify the corresponding cell datum as in Definition \ref{definition gca}. Take
 $$ \Lambda = \Lambda (\sub{w}) := \{ x \in W \mbox{ } | \mbox{ } w \geq x \}$$ partially ordered by reversing the usual Bruhat order. Accordingly, $w$ and $e$ (where $e$ denotes the identity of $W$) are the minimal and the maximal element in $\Lambda (\sub{w})$, respectively. For each $x \in \Lambda(\sub{w})$, define $T(x) := \mathbb{L}_{\underline{w}}(x)$, i.e., $T(x)$ is the set of all leaves in $\mathbb{T}_{\sub{w}}$ with final target $x$. We end the description of the cell datum by defining $c_{l_1l_2}^{x} =: l_1^{a}\circ l_2$, for all $l_1, l_2 \in T(x)$ and $x\in \Lambda$.

\medskip
On the other hand, there is a natural degree function
\begin{equation}   \label{degree function}
\deg: \coprod_{x\in \Lambda(\sub{w})} T(x) \rightarrow \Z
\end{equation}
given by the degree of the leaves. The reader might expect that the double leaves basis is a graded cellular basis for $\en^{\Z}(\bs{w})$. However, this is not true because $\en^{\Z}(\bs{w})$ is not a graded algebra in the sense of the definition in Section \ref{section gca}. Actually, $\en^{\Z}(\bs{w})$ does not satisfy the property $A_iA_j \subset A_{i + j}$, for all $i, j \in \Z$, since the ground ring $R$ is not of degree zero. This drawback can be rectified by reducing modulo $R^{+}$. We recall that $\R \simeq R/R^{+}$ and  define $\A = \en^{\Z}(\bs{w})\otimes_{R} \mathbb{R} $.

\begin{teo}  \label{graded celullarity}
 The set $\{l_1^{a} \circ l_{2}\otimes_{R} 1 \mbox{ } | \mbox{ } l_{1}, l_{2} \in \mathbb{L}_{\sub{w}}(x); x \in \Lambda (\sub{w})\}$ is a graded cellular basis for $\A$.
\end{teo}
\begin{dem}
The cellularity of $\A$ is clear from the cellularity of $\en^{\Z}(\bs{w})$. For the graded part, we only need to check that
$$ \deg (l_1^{a} \circ l_{2}) = \deg (l_{1}) + \deg(l_2)$$
which follows directly from the definition of a double leaf.
\end{dem}

\medskip
Given the details of the cellular structure of $\A$ we have automatically  defined the corresponding graded cell $\A$-modules and graded simple $A$-modules, as well as the graded  decomposition numbers for $\A$. However, by the abstract definition of cell modules and their bilinear forms given in Section \ref{section gca}, it is not clear how one ought to work with them. Fortunately, in this case, we can be a little more specific. Let us provide two definitions.

\begin{defi} \label{factors through} \rm
Let $\sub{w}$ and $\sub{v}$ be expressions and $u \in W$. We say that a double leaf $l_1^a \circ l_2 \in \mathbb{L}_{\sub{v}}^a  \cdot \mathbb{L}_{\sub{w}}$ \emph{factors through} $u$ if $l_1 \in \mathbb{L}_{\sub{v}}(u)$ and $l_2 \in\mathbb{L}_{\sub{w}}(u) $.
\end{defi}

\begin{defi} \label{ideal} \rm
Let $\sub{w}$ and $\sub{v}$ be expressions. For $u \in W$, we define the set $\mathbb{DL}_{ < u}$ as the spans of the double leaves in $\mathbb{L}_{\sub{v}}^a  \cdot \mathbb{L}_{\sub{w}}$ that factor through $x < u$.
\end{defi}

Let us denote by $\cell{x}$, $\simple{x}$, and $\gdn{x}{y}$ the graded cell modules, the graded simple modules, and the graded decomposition numbers of $\A$, for $x, y \in \Lambda(\sub{w})$, respectively, corresponding to the cellular structure determined by the double leaves basis. We now explicate the action of $\A$ on a cell module. Let $x \in \Lambda (\sub{w})$. By definition, the graded cell module $\cell{x}$ is the $\R$-vector space spanned by $\mathbb{L}_{\sub{w}}(x)$.

\begin{rem} \rm
If we want to be completely consistent with the notation introduced in Section \ref{section gca}, the cell module must be the $\R$-vector space with basis $$ \{ c_{l}^{x} \mbox{ } | \mbox{ } l \in  \mathbb{L}_{\sub{w}}(x) \}. $$ However, to avoid a subindex catastrophe, we prefer the previous notation.
\end{rem}

Let $a \in\en^{\Z}(\bs{w})$ and $l \in \mathbb{L}_{\sub{w}}(x)$. To determine $l(a\otimes 1) \in \cell{x}$, we calculate the expansion of $l\circ a$ in terms of the double leaves basis for $ \hm^{\Z}(\bs{w}, \bs{x} ) $. It is not difficult to note that
\begin{equation}
l\circ a \equiv \sum_{g \in\mathbb{L}_{\sub{w}}(x) }  gr_{g}  \mod \mathbb{DL}_{<x}
\end{equation}
for some scalars $r_{g}\in R$. Then, the action of $\A$ on $\cell{x}$ is
\begin{equation}
l( a\otimes 1) = \sum_{g\in\mathbb{L}_{\sub{w}}(x) }  g(r_{g}\otimes 1) \in \cell{x}.
\end{equation}
In a similar manner, we can explicate the bilinear form on $\cell{x}$ induced by the cellular structure. Let $l_1, l_2 \in \mathbb{L}_{\sub{w}}(x)$. Now, $l_1\circ l_2^{a}\in \en^{\Z}(\bs{x})$. Thus, we can expand it in terms of the double leaves basis for $\en^{\Z}(\bs{x})$. Again, it is not hard to note that
\begin{equation}    \label{description bilinear form}
 l_1^{a}\circ l_2\equiv \mathbb{I}_{x} r(l_1, l_2)  \mod \mathbb{DL}_{<x}
\end{equation}
for some  $r(l_1, l_2) \in R $, and  where $\mathbb{I}_{\sub{x}}$ denotes the identity map of $\bs{x}$. Then, the value of the bilinear form $\lr{ \text{ , }}$ on $\cell{x}$ at two leaves $l_1$ and $l_2$ is $\lr{l_1, l_2} = r(l_1, l_2) \otimes 1 \in \R$. Since $\deg (\mathbb{I}_{x}) = 0$, we have
 \begin{equation}
 \deg (l_1) + \deg (l_2) =  \deg (r(l_1, l_2)).
\end{equation}
Thus, $\lr{l_1, l_2} = r(l_1, l_2) \otimes 1 = 0 $ unless $\deg (l_1) + \deg (l_2) = 0$.

\medskip
It is a straightforward exercise to confirm that the descriptions of cell modules and bilinear form provided here coincide with those in Section \ref{section gca}. Let us denote by $\Lambda_{0}(\sub{w})$ the set that parameterizes the entire set (up to degree shift) of simple modules of $\A$, i.e.,
\begin{equation}
\Lambda_{0}(\sub{w}) = \{ x \in \Lambda(\sub{w}) \mbox{ }| \mbox{ } \simple{x}  \neq 0 \}
\end{equation}
In order to obtain an interpretation of the KL-polynomials as graded decomposition numbers, we need the following two lemmas:

\begin{lem}  \label{graded dimensions}
Let $\sub{w}$ be a reduced expression of $w \in W$. Then, for all $x \in \Lambda(\sub{w})$
\begin{equation} \label{graded dimension cell}
\gd{\cell{x}} = \sum_{l \in \mathbb{L}_{\sub{w}}(x)} v^{\deg(l)}= \gd \hm^{\Z}(\bs{w}, R_{x} )\otimes_{R} \R
\end{equation}
Furthermore, the coefficient of $v^{k}$ in $\gd{\simple{y}}$ is the multiplicity $B_{y} \langle k \rangle$ as a direct summand in $\bs{w}$, for all $y \in \Lambda(\sub{w})$. We thus have
\begin{equation}
\bs{w} \cong  \bigoplus_{y \in \simpleset{w}} \gd{D_{\sub{w}}(y)} B_{y}
\end{equation}
\end{lem}

\begin{dem}
The left side of (\ref{graded dimension cell}) is clear from the definitions, and the right side follows from Lemma \ref{standard basis}. The last claim is a direct consequence of \cite[Lemma 3.1]{wi}, and the description of the bilinear form given in (\ref{description bilinear form}) and its homogeneity.
\end{dem}

\begin{lem}   \label{isomorphism cell}
Let $\sub{w}$ be a reduced expression for $w \in W$. Then, for all $x \leq w$, we have an isomorphism
\begin{equation}  \label{claim}
\hm^{\Z}(\bs{w}, R_x)\otimes_{R} \R \cong \cell{x}
\end{equation}
of right $\A$-modules.
\end{lem}
\begin{dem}
Note first that $\hm^{\Z}(\bs{w}, R_x)\otimes_{R} \R$ has a natural structure of a right $\A$-module by composition of morphisms. Concretely, if $g \in \hm^{\Z}(\bs{w}, R_x)$  and $a \in \en^{\Z}(\bs{w})$, the action of $\A$ on $\hm^{\Z}(\bs{w}, R_x)\otimes_{R} \R$ is given by
\begin{equation}
(g \otimes 1) (a \otimes 1) = (g \circ a) \otimes 1
\end{equation}
Further, by Lemma \ref{standard basis}, $\hm^{\Z}(\bs{w}, R_x)\otimes_{R} \R$ is an $\R$-vector space with basis
$$\{ (\beta_x \circ l )\otimes 1 \mbox{ } | \mbox{ }l \in  \mathbb{L}_{\sub{w}} (x) \},$$
where $\beta_x: \bs{x} \rightarrow R_x$ is the bimodule morphism defined following Lemma \ref{lemma lex}. Since $\cell{x}$ is defined as the $\R$-vector space with basis $\mathbb{L}_{\sub{w}}(x)$,  there is a canonical $\R$-linear isomorphism determined by
\begin{equation}
\begin{array}{cccc}
  f: & \cell{x} & \rightarrow & \hm^{\Z}(\bs{w}, R_x)\otimes_{R} \R  \\
   & l &  \rightarrow &   (\beta_x \circ l )\otimes 1
\end{array}
\end{equation}
for all $l\in \mathbb{L}_{\sub{w}}(x) $. Then, by the $\R$-linearity of $f$, to finish the proof we need to show that
\begin{equation}  \label{equation f}
f(l (a\otimes 1))=f(l)(a\otimes 1),
\end{equation}
 for all $a\in \en^{\Z}(\bs{w})$ and $l\in\mathbb{L}_{\sub{w}}(x) $. We prove that both sides of (\ref{equation f}) are equal to $(\beta_{x}\circ l \circ a) \otimes 1$. First, note that
  $$ f(l)(a\otimes 1) = ( (\beta_x \circ l )\otimes 1 ) (a\otimes 1)=(\beta_{x}\circ l \circ a) \otimes 1,  $$
proving that the right side of (\ref{equation f}) is equal to $(\beta_{x}\circ l \circ a) \otimes 1$. To prove the other equality, we need the following

\begin{cla}  \label{claim beta}  \rm
If $g \in \hm^{\Z}(\bs{w}, \bs{x})$ belongs to $\mathbb{DL}_{< x}$ then $\beta_{x} \circ g = 0$.
\end{cla}
\begin{dem}
It is enough to show that
\begin{equation}  \label{beta kills}
\beta_{x} \circ (l_2^a \circ l_1) =0,
\end{equation}
for all double leaves $(l_2^a \circ l_1) \in \hm^{\Z}(\bs{w}, \bs{x})$ that factor through $ u < x $. On the contrary, suppose that there exists a double leaf $(l_2^a \circ l_1) \in \hm^{\Z}(\bs{w}, \bs{x})$ that factors through $ u < x $ such that $\beta_{x} \circ (l_2^a \circ l_1) \neq 0$. In particular, we have $\beta_{x} \circ l_2^a \neq 0 $. Note that $\beta_{x} \circ l_2^a$ belongs to $ \hm^{\Z}(\bs{u}, R_x)$, for some reduced expression $\sub{u}$ of u. Therefore,  $ \hm^{\Z}(\bs{u}, R_x) \neq 0 $. This contradicts  Corollary \ref{corollary standard basis} since $u < x$, proving (\ref{beta kills})and  Claim \ref{claim beta}.
\end{dem}

\medskip
Let us return to the proof of the lemma. To conclude the proof, we need to show $f(l (a\otimes 1)) = (\beta_{x}\circ l \circ a) \otimes 1$. Write
\begin{equation}   \label{eq a}
l \circ a \equiv \sum_{g \in \mathbb{L}_{\sub{w}}(x)} gr_{g} \mod \mathbb{DL}_{<x},
\end{equation}
for some scalars $r_{g} \in R$. Composing with $\beta_{x}$ to the left in the previous equation and using Claim \ref{claim beta}, we obtain
\begin{equation}   \label{eq b}
\beta_x \circ l \circ a = \sum_{g \in \mathbb{L}_{\sub{w}}(x)} (\beta_{x} \circ g)r_{g}.
\end{equation}
Thus, by reducing modulo $R^{+}$ we obtain
\begin{equation}    \label{eq c}
(\beta_x \circ l \circ a)\otimes 1  = \sum_{g \in \mathbb{L}_{\sub{w}}(x)} (\beta_{x} \circ g)r_{g}\otimes 1 \in \hm^{\Z}(\bs{w}, R_x)\otimes_{R} \R.
\end{equation}
On the other hand, by (\ref{eq a}) we know
\begin{equation}    \label{eq d}
l (a\otimes 1) = \sum_{g \in \mathbb{L}_{\sub{w}}(x)}  g(r_{g} \otimes 1) \in \cell{x}.
\end{equation}
Thus,
\begin{equation}      \label{eq e}
f(l (a\otimes 1)) = \sum_{g \in \mathbb{L}_{\sub{w}}(x)}  (\beta_{x} \circ g)(r_{g} \otimes 1) \in \cell{x}.
\end{equation}
Combining (\ref{eq c}) with (\ref{eq e}), we conclude that $f(l (a\otimes 1)) = (\beta_{x}\circ l \circ a) \otimes 1$. This completes the proof of the lemma.
\end{dem}

\medskip
We are now in a position to interpret the Kazhdan-Lusztig polynomials as graded decomposition numbers. This result is the key to proving the monotonicity conjecture for coefficients of the Kazhdan-Lusztig polynomials in the following section.

\begin{teo}  \label{teo hl as gdn}
 Let $\sub{w}$ be a reduced expression for $w \in W$ and $x \leq w$. Then,
 \begin{equation}   \label{eq kl as gdn}
d_{\sub{w}}(x, w) = h_{x, w}
 \end{equation}
\end{teo}
\begin{dem}
By Lemma \ref{graded dimensions}, we have the following isomorphism
\begin{equation} \label{bs decomposition}
\bs{w} \cong  \bigoplus_{y \in \simpleset{w}} \gd{D_{\sub{w}}(y)} B_{y}
\end{equation}
as $(R, R)$-bimodules. Since $\mathbb{L}_{\sub{w}}(w) = \{ \mathbb{I}_{\sub{w}} \}$, it is easy to note that $ \simple{w} = \cell{w}$,
where $ \mathbb{I}_{\sub{w}}$ is the identity map on $\bs{w}$. Therefore, $w \in \simpleset{w}$ and $\gd{D_{\sub{w}}(w)} = 1$. Now, we can choose a projector (idempotent) $e \in \en^{\Z}( \bs{w})$ whose image is isomorphic to $B_{w}$. We denote by $\hat{e} \in \A$ its reduction modulo $R^{+}$. Note that $e$ and $\hat{e}$ are primitive idempotents. By (\ref{bs decomposition}), we have the following isomorphism of right $R$-modules
\begin{equation} \label{hom bs decomposition}
\hm^{\Z}(\bs{w}, R_x) = \bigoplus_{y \in \simpleset{w}} \gd{\simple{y}}\hm^{\Z} (B_{y}, R_{x})
\end{equation}
for all $x \in W$. Composing using $e$ from the right, we obtain
$$\hm^{\Z}(\bs{w}, R_x) e \simeq \hm^{\Z} (B_{w}, R_{x})$$
as right $R$-modules. Hence,
$$(\hm^{\Z}(\bs{w}, R_x)\otimes_R  \R )\hat{e} \simeq \hm^{\Z} (B_{w}, R_{x}) \otimes_R \R $$
as graded $\R$-vector spaces. Taking the graded dimension on both sides, we obtain
\begin{equation}
\gd{ \hm^{\Z}(\bs{w}, R_x)\otimes_R  \R )\hat{e}} = \gd{ \hm^{\Z} (B_{w}, R_{x}) \otimes_R \R } = h_{x, w}
\end{equation}
where the last equation is from Soergel's conjecture and Theorem \ref{inverse categorification}. Therefore, by Lemma \ref{isomorphism cell}, we have
\begin{equation} \label{equation kl}
\gd \cell{x} \hat{e} =  h_{x, w}.
\end{equation}
On the other hand, by Theorem \ref{idempotent algebra}, the algebra $\hat{e} \A \hat{e}$ is a graded algebra with a unique graded simple module (up to degree shift). Actually, $\simple{w}\hat{e}$ is the unique simple $\hat{e} \A \hat{e}$-module since
\begin{equation}  \label{graded dimension equal to one}
 \gd \simple{w}\hat{e} = \gd \cell{w} \hat{e} = h_{w, w} = 1.
\end{equation}
Finally, we obtain
$$
\begin{array}{rl}
  h_{x, w}         &  = \gd \cell{x} \hat{e}     \\
         &  =  d(\cell{x} \hat{e}, \simple{w} \hat{e})      \\
         &  =  d(\cell{x},\simple{w})  \\
         &  = \gdn{x}{w}
\end{array}
$$
where the second equation follows from (\ref{graded dimension equal to one}) and Theorem \ref{idempotent algebra}(d), and the third  equation follows from Theorem \ref{idempotent algebra}(b).
\end{dem}

\medskip
Note that the left side of (\ref{eq kl as gdn}) depends on the expression $\sub{w}$ whereas the right side does not. Thus, the theorem also claims that $\gdn{x}{w}$ does not depend on the choice of the reduced expression $\sub{w}$ of $w$. Note also that $\gdn{x}{w}$ is only defined for $x\leq w$ whereas the KL-polynomials are defined for each pair of elements in $W$.  However, the above is irrelevant because the KL-polynomial $h_{x, w} \neq 0 $ if and only if $x \leq w$. Summing up, the above theorem says that each nonzero KL-polynomial can be interpreted as a graded decomposition number.

\section{Monotonicity}  \label{section monotony}
In this section, we prove the  Monotonicity Conjecture for the coefficients of the Kazhdan-Lusztig polynomials. More precisely, we prove:
\begin{con}  \label{conj}
Let $W$ be any Coxeter group. If $u, v, w \in W$ and $u \leq v \leq w$ then
\begin{equation}
P_{u, w}(q) - P_{v, w}(q) \in \mathbb{N}[q]
\end{equation}
\end{con}

In terms of the polynomials $h_{x, w}(v) \in \mathbb{Z}[v]$, the above conjecture is equivalent (via Remark \ref{variables}) to
\begin{equation}  \label{eq with h}
h_{u, w}(v) - v^{l(v) - l(u)}h_{v, w}(v) \in \mathbb{N}[v]
\end{equation}

We prove (\ref{eq with h}) in this section. To do this, we are first interested in a particular leaf.


\begin{lem}  \label{largest leaf}
Let $W$ be a Coxeter group. Let $u, v \in W$ with $u \leq v$ and let $\sub{v} $ be a reduced expression for $v$. Then, there is a unique leaf in $\mathbb{L}_{\sub{v}}(u)$ of degree $l(v) - l(u)$.
\end{lem}

\begin{dem}
This is a direct consequence of the definition of the leaves and \cite[Lemma 5.1]{li} or \cite[Proposition 2.3]{deo}.
\end{dem}

\medskip
Let $u, v, w \in W$ with $u \leq v \leq w$. For the rest of the paper, we fix reduced expressions $\sub{u}$, $\sub{v}$, and $\sub{w}$  for $u$, $v$, and $w$, respectively. We denote by $G_{v}^{u}$ the leaf in Lemma \ref{largest leaf}, and refer to it as the \emph{largest leaf} from $\bs{v}$ to $\bs{u}$. It follows directly from the construction of the leaves that for all $l \in \mathbb{L}_{\sub{v}}(u)$,
\begin{equation}   \label{largest leaf eq one}
\begin{array}{rl}
  \deg (l) = & n_{m}(l) - n_{j}(l)  \\
 l(v) - l(u) =  &  n_{m}(l) - n_{j}(l),
\end{array}
\end{equation}
where $n_{m}(l)$ (resp. $n_{j}(l)$) denotes the number of times that morphisms of type $m_s$ (resp. $j_s$) appear in the construction of leaf $l$. In particular, if we set $l = G_{v}^{u}$ in (\ref{largest leaf eq one}) and subtract the resulting equations, we obtain
\begin{equation}  \label{largest leaf eq three}
n_{j}(G_{v}^u) = 0,
\end{equation}
since $\deg ( G_{v}^{u} ) = l(v) - l(u)$. That is, morphisms of type $j_s$ do not appear in the construction of the largest leaf. We define a map, $\Phi_{\sub{w}}^{u, v}$, from $\hm^{\Z}(\bs{w}, R_v) \otimes_{R} \R$ to $\hm^{\Z}(\bs{w}, R_u) \otimes_{R} \R$, determined in the standard basis by


\begin{equation}
\begin{array}{cccc}
  \Phi_{\sub{w}}^{u, v}: & \hm^{\Z}(\bs{w}, R_v) \otimes_{R} \R  & \rightarrow & \hm^{\Z}(\bs{w}, R_u) \otimes_{R} \R  \\
& (\beta_{v} \circ l )\otimes 1 & \rightarrow &  (\beta_{u} \circ G_{v}^{u} \circ l ) \otimes 1
\end{array}
\end{equation}
for all $l \in \mathbb{L}_{\sub{w}}(v)$. The map $\Phi_{\sub{w}}^{u, v}$ is the key to demonstrating the monotonicity conjecture at the conclusion of this section. In order to know the properties of  $\Phi_{\sub{w}}^{u, v}$, we need a notational and technical lemma.

\begin{defi}  \rm
Let $\sub{w}$ be an expression. For $b \in \bs{w}$, we define $\co{b}$ as the coefficient of $1^{ \otimes}$ in the expansion of $b$ in terms of the basis of $\bs{w}$ described in (\ref{base of bott-samelson}).
\end{defi}

\begin{lem}   \label{technical lemma}
Let $u, v, w \in W$ with $u \leq v \leq w$. If $h \in \hm^{\Z}(\bs{w}, \bs{v})$ belongs to  $\mathbb{DL}_{< v}$, then,
\begin{equation}
G_{v}^{u} \circ h \equiv \sum_{g\in \mathbb{L}_{\sub{w}}(u)}  gr_{g} \mod \mathbb{DL}_{<u},
\end{equation}
for some scalars $r_{g} \in R^{+}$.
\end{lem}
\begin{dem}
Let $l_2^{a}\circ l_{1} $ be a double leaf morphism in $\hm^{\Z}(\bs{w}, \bs{v})$ that factors through $z < v$. Write
\begin{equation}  \label{expansion}
 G_{v}^{u} \circ (l_2^{a}\circ l_{1}) = \sum_{g \in \mathbb{L}_{\sub{w}}(u)}  gr_{g}  +   f ,
\end{equation}
for some scalars $r_{g} \in R$ and some morphism $f \in \mathbb{DL}_{< u}$. To finish the proof, we need to show that $r_{g} \in R^{+}$, for all $g \in \mathbb{L}_{\sub{w}}(u)$. Recall the indexation (given in Section \ref{section dlb}) of the leaves by two sequences of zeros and ones, and define
\begin{equation}
\mathcal{J} = \{ \bol{j} \in \{0, 1\}^{l(\sub{w})} \mbox{ } |\text{ there is a leaf } g \in \mathbb{L}_{\sub{w}}(u) \text{ such that } g = f_{\bol{i}}^{\bol{j}} \}
\end{equation}
By \cite[Lemma 5.1]{li}, we know that each $\bol{j} \in \mathcal{J}$ determines a unique leaf in $\mathbb{L}_{\sub{w}}(u)$. Index $\mathcal{J} = \{\bol{j}_{1}, \ldots , \bol{j}_{n}  \}$ so that $\bol{j}_{k} < \bol{j}_{m}$ ($<$ here denotes the lexicographical order) if and only $k < m$, for all $1 \leq k < m \leq n$. For $1 \leq k \leq n$, we denote by $g_{k}$ the leaf determined by $\bol{j}_{k}$. On the other hand, by the construction of the leaves and (\ref{largest leaf eq three}), it is easy to note that
\begin{equation} \label{they are in R plus}
\co{ (G_{v}^{u} \circ l_2^{a}\circ l_{1})(b)} \text{, } \co{f(b)}  \in R^{+},
\end{equation}
for all $b \in \bs{w}$. We now proceed by induction. If we evaluate (\ref{expansion}) at $x^{\bol{j}_1}$, then by Lemma \ref{lemma lex}  and (\ref{they are in R plus}), we find that
\begin{equation} \label{they are in R plus  one}
r_{g_{1}} =   \co{ (G_{v}^{u} \circ l_2^{a}\circ l_{1})(x^{\bol{j}_1})} - \co{f(x^{\bol{j}_1})}  \in R^{+},
\end{equation}
which provides the basis of our induction. Now, let $1 < k \leq n$ and assume that we have already proved that $r_{g_{m}} \in R^{+}$, for all $1 \leq m <k$. Evaluating (\ref{expansion}) at $x^{\bol{j}_k}$, and again using \cite[Lemma 5.1]{li}, we obtain
\begin{equation} \label{they are in R plus  two}
r_{g_{k}} =   \co{ (G_{v}^{u} \circ l_2^{a}\circ l_{1})(x^{\bol{j}_k})} - \sum_{m = 1}^{k - 1}  \co{ g_{m}(x^{\bol{j}_k})  } r_{g_{m}} - \co{f(x^{\bol{j}_k})}.
\end{equation}
By (\ref{they are in R plus}) and the inductive hypothesis, we know that the right side of (\ref{they are in R plus  two}) belongs to $R^{+}$. Therefore, $r_{g_{k}} \in R^{+}$. This completes the induction and the proof of the lemma.
\end{dem}

\begin{pro}   \label{Phi is hom}
Let $u, v, w \in W$ with $u \leq v \leq w$. Then, $ \Phi_{\sub{w}}^{u, v}$ is a homogeneous $\A$-module homomorphism of degree $l(v) - l(u)$.
\end{pro}
\begin{dem}
The claim that $ \Phi_{\sub{w}}^{u, v}$ is homogeneous with degree  $l(v) - l(u)$ is a direct consequence of the definitions as well as the fact that $\deg (G_{v}^{u}) = l(v) - l(u)$. Now, in order to prove that $ \Phi_{\sub{w}}^{u, v}$ is an $\A$-module homomorphism, it is sufficient to show that
\begin{equation}  \label{eq Phi one}
\Phi_{\sub{w}}^{u, v}(((\beta_{v} \circ l) \otimes 1)(a \otimes 1) ) = \Phi_{\sub{w}}^{u, v}((\beta_{v} \circ l) \otimes 1 )(a \otimes 1)
\end{equation}
for all $l \in \mathbb{L}_{\sub{w}}(v)$ and $a \in \en^{\Z}(\bs{w}) $. We prove that both sides of (\ref{eq Phi one}) are equal to $(\beta_{u}\circ G_{v}^{u} \circ l \circ a) \otimes 1$. The desired equality for the right side of (\ref{eq Phi one}) is trivial because by the definition of $\Phi_{\sub{w}}^{u, v}$, we have
\begin{equation}   \label{eq Phi two}
\begin{array}{rl}
  \Phi_{\sub{w}}^{u, v}((\beta_{v} \circ l) \otimes 1 )(a \otimes 1) & =  ((\beta_{u} \circ G_{v}^{u} \circ l ) \otimes 1)(a \otimes 1) \\
   & = (\beta_{u} \circ G_{v}^{u} \circ l \circ a ) \otimes 1
\end{array}
\end{equation}
To obtain the equality for the left side of (\ref{eq Phi one}), we first write
\begin{equation}   \label{eq que faltaba}
l\circ a \equiv \sum_{f\in \mathbb{L}_{\sub{w}}(v)} fr_{f} \mod \mathbb{DL}_{<v},
\end{equation}
for some scalars $r_{f} \in R$.  By Claim \ref{claim beta}, we know that
\begin{equation}  \label{falta equation}
((\beta_{v} \circ l) \otimes 1)(a \otimes 1) = \sum_{f\in \mathbb{L}_{\sub{w}}(v)} \beta_{v}\circ fr_{f} \otimes 1
\end{equation}
Thus, by applying $\Phi_{\sub{w}}^{u, v}$ to \ref{falta equation}, we have
\begin{equation}  \label{eq Phi three}
\Phi_{\sub{w}}^{u, v}(((\beta_{v} \circ l) \otimes 1)(a \otimes 1) ) =    \sum_{f\in \mathbb{L}(v)} (\beta_{u} \circ G_{u}^{v} \circ fr_{f} ) \otimes 1.
\end{equation}
On the other hand, by composing (\ref{eq que faltaba}) with $G_{u}^{v}$ to the left and using Lemma \ref{technical lemma}, we obtain
\begin{equation}  \label{eq Phi four}
G_{v}^{u} \circ l\circ a \equiv \sum_{f\in \mathbb{L}_{\sub{w}}(v)} G_{v}^{u} \circ fr_{f} + \sum_{g\in \mathbb{L}_{\sub{w}}(u)}g\rho_{g}   \mod \mathbb{DL}_{<u},
\end{equation}
for some scalars $\rho_{g}\in R^{+}$. Now, by composing with $\beta_{u}$ to the left in (\ref{eq Phi four}) and using Claim \ref{claim beta}, we have
\begin{equation}  \label{eq Phi five}
\beta_{u} \circ G_{v}^{u} \circ l\circ a = \sum_{f\in \mathbb{L}_{\sub{w}}(v)}\beta_{u} \circ G_{v}^{u} \circ  fr_{f} + \sum_{g\in \mathbb{L}_{\sub{w}}(u)} \beta_{u} \circ g\rho_{g}   .
\end{equation}
Following this, by reducing modulo $R^{+}$ and using the fact that $\rho_{g} \in R^{+}$ for all $g\in \mathbb{L}_{\sub{w}}(u)$, we obtain
\begin{equation}  \label{eq Phi six}
\beta_{u} \circ G_{v}^{u} \circ l\circ a  \otimes 1 = \sum_{f\in \mathbb{L}_{\sub{w}}(v)}\beta_{u} \circ G_{v}^{u} \circ  fr_{f} \otimes 1.
\end{equation}
Finally, combining (\ref{eq Phi three}) with (\ref{eq Phi six}), we obtain
$$\Phi_{\sub{w}}^{u, v}(((\beta_{v} \circ l) \otimes 1)(a \otimes 1) ) =  \beta_{u} \circ G_{v}^{u} \circ l\circ a  \otimes 1. $$
This completes the proof of the proposition.
\end{dem}

\begin{pro}  \label{injective hom}
 Let $u, v, w \in W$ with $u \leq v \leq w$. The map $\Phi_{\sub{w}}^{u, v}$ is injective.
\end{pro}

\begin{dem}
Let us suppose that
\begin{equation}  \label{phi injec}
\Phi_{\sub{w}}^{u, v}\left(  \sum_{l\in \mathbb{L}_{\sub{w}}(v)} (\beta_{v} \circ lr_{l}) \otimes 1 \right) = 0,
\end{equation}
for some $r_{l} \in R$. The proof is completed by showing that $r_{l}\in R^{+}$, for all $l\in \mathbb{L}_{\sub{w}}(v)$. By the definition of $\Phi_{\sub{w}}^{u, v}$, Equation (\ref{phi injec}) implies
\begin{equation}  \label{phi injec two}
 \sum_{l\in \mathbb{L}_{\sub{w}}(v)} (\beta_{u} \circ  G_{v}^{u} \circ lr_{l}) \otimes 1  = 0,
\end{equation}
Hence
\begin{equation}  \label{phi injec three}
 \sum_{l\in \mathbb{L}_{\sub{w}}(v)} (\beta_{u} \circ  G_{v}^{u} \circ lr_{l})(b) \in R^{+},
\end{equation}
for all $b \in \bs{w}$. Now, we can prove that $r_{l}\in R^{+}$, for all $l \in \mathbb{L}_{\sub{w}}(v)$ by the same inductive method as used in the proof of Lemma \ref{technical lemma}. In this context, (\ref{phi injec three}) takes the place of (\ref{they are in R plus}). The details are left to the reader.
\end{dem}

\begin{teo}
Conjecture \ref{conj} holds.
\end{teo}
\begin{dem}
A direct consequence of Lemma \ref{isomorphism cell}, Proposition \ref{Phi is hom} and Proposition \ref{injective hom} is
 that $\cell{v} \langle l(v) - l(u) \rangle $ is a graded right $\A$-submodule of $\cell{u}$. Therefore,
\begin{equation} \label{ultima}
\gdn{u}{x} - v^{l(v) - l(u)}\gdn{v}{x} \in \mathbb{N}[v],
\end{equation}
for all $x \in \simpleset{w}$. In particular, since we know that $w \in \simpleset{w}$, (\ref{ultima}) holds for $x = w$. By Theorem  \ref{teo hl as gdn}, for $x = w$, (\ref{ultima}) becomes (\ref{eq with h}), thus proving the theorem.
\end{dem}

\section*{acknowledgments}
This research was supported by a Fondecyt Postdoctoral grant no. $3140612$. The author would like to thank Nicolas Libedinsky for many stimulating conversations on topics concerning this paper.

\begin{scriptsize}

\end{scriptsize}

\end{document}